\documentclass{article}
\usepackage{amsbsy,amssymb,amscd,amsfonts,latexsym,amstext,delarray,
amsmath}  \headsep=15pt
\def\A{ {\bf A} }

\def\ra{{ \rightarrow }}

\def\C{{ \mathbb{C} }}

\def\Z{{ \mathbb{Z}}}

\newtheorem{thm}{Theorem}[section]

\newtheorem{conj}[thm]{Conjecture}

\newtheorem{observation}[thm]{Observation}

\newtheorem{guess}[thm]{Guess}
\usepackage{amsmath}
\usepackage{amsfonts}
\usepackage{amscd}
\usepackage{amssymb}
\def\Q{\mathbb{Q}}

\def\P{ {\bf P}}

\title{On relative computability for curves}
\author{ Minhyong Kim}

\begin{document}
\maketitle

This semester, I've enjoyed greatly the hospitality of Carl Jockusch who was gracious enough to share his
office with me for the duration of my stay at UIUC. About two weeks ago, I also had the pleasure of conversing
with him over lunch. It was then that I  learned of  the term `relative computability,' a subject to which I
am told Carl has made profound contributions. It was thus, with all the naive enthusiasm of an amateur, that I
thought to include it in the title of my lecture. Furthermore, Carl was kind enough to inform me of a
remarkable conjecture  that provides a comfortable framework for the topic I am set to discuss. So let us
start with that. Here, the objects of interest are integral Diophantine equations, that is, equations of the
form
$$f(x_1, x_2, \ldots, x_n)=0$$
where $f(x_1,\ldots, x_n)$ is a polynomial with integral coefficients. I'm sure you are well aware of the {\em
existence problem} related to this equation, that is, the problem of existence of integral solutions, as well
the undecidability result of Matiyasevich. The subject of the conjecture, however, goes beyond existence. That
is, it considers simultaneously the problem of determining the {\em finiteness} (at most) of the solution set.
It is easy to see, by adjoining  a dummy variable, that the finiteness
problem is undecidable as well. However, as you
know, some of the most celebrated theorems of arithmetic are about finiteness
for specific sorts of equations. In fact, many of them state
finiteness  in total ignorance of existence. And then, sometimes you know
existence and nothing about finiteness.
But as far as the decision problem is concerned,
the conjecture in question probes this relationship more
deeply.
\begin{conj}
The finiteness problem for integral points
is undecidable relative to the existence problem.
\end{conj}
In other words, even given an `existence oracle', i.e., a decision oracle for the existence problem (or
equivalently, an oracle for the halting problem), the finiteness problem should be undecidable. I am sure one
can make this conjecture more precise or generalize it in many ways using the sophisticated machinery of
recursion theory, of which I am woefully ignorant. Also, for a naive number-theorist, the subtleties of
relative computability are often hard to comprehend in a situation where the oracle whose existence we need to
assume is known not to exist. This is of course because we are more obsessed  with solving problems than
classifying them. In any case, when I heard this conjecture, it seemed natural enough to ask this question in
a context where a relative computability result still has a chance of leading to an actual reduction of the
problems of interest. That is, we can shift our attention (seemingly slightly) to {\em rational solutions}
rather than just the integral ones. Having done that, one finds that this conjecture relates rather well to
established programs in Diophantine geometry, and a precise articulation of this relationship becomes quite
desirable. I will not attempt to carry this out today, out of pure laziness. However, I do wish to give some
sense of the issues that come up, and maybe put forth a suggestion or two as to the kind of phenomena one
should expect. For example, here is something of a
\begin{guess}
For rational solutions, the finiteness problem is {\em decidable} relative to the existence problem.
\end{guess}
That is, as far as rational solutions are concerned, my expectation goes counter to the conjecture for integral
solutions. That this {\em could} be so  is not at all surprising  since, I believe, experienced recursion
theorists and number theorists do find the nature of rational and integral solutions to be very different. In
particular, I don't think too many people expect $\Z$ to be definable in $\Q$. Barry Mazur has pointed out
that the explicit family of equations constructed by Matiyasevich all have rational solutions for trivial
reasons (I have not verified this myself). As another illustration, consider the case of elliptic curves (curves of genus one equipped with a rational point) where the
finiteness of integral solutions is well-known. For rational solutions, by contrast, the decision
problem for finiteness forms an important part of the Birch and Swinnerton-Dyer (BSD) conjecture. In other
words, the rational case is notoriously difficult. More on this point later.

I wish to postpone to the end of the lecture an explanation of the intuition behind my guess which needs to be
vague anyways, since otherwise, I would have done the work necessary to elevate the guess to a conjecture.
Instead, I will concentrate for now on the case of curves. That is, to say, we are interested in Diophantine
equations in  two-variables,
$$f(x,y)=0,$$
where $f$ is again a polynomial with integral coefficients, but assumed now to be irreducible over the complex
numbers.  [A note for geometers: In this lecture, although I can't help lapsing into geometric terminology,
the emphasis really is on the equations themselves. That is, the precise presentation under discussion, as
input for  machines, is the focal point.] As stated, we will be interested in the rational
solutions, which I will mostly refer to merely as solutions, for brevity. A rough classification of the
solution set, representing the main achievements of 20th century number theory, depends on the {\em genus},
$g(f)$, of the equation (or the polynomial). That is, one considers the field $\C(x)[y]/(f(x,y))$ which can be
realized as the field of meromorphic functions on a unique Riemann surface. The genus of this surface is the
genus of $f$. In most cases, it can be computed readily from $f$ using the formula $$g(f)=(d-1)(d-2)/2$$ where
$d$ is the degree of $f$. And then one knows:

If $g(f)=0$ then the solution set is empty or infinite.

If $g(f)=1$ then the solution set can be empty, non-empty finite, or infinite. In short, anything is possible.

Finally, if $g(f)\geq 2$, then the solution set is empty, or non-empty finite.

A  fact that emerges from this classification is that if we restrict our attention to equations with $g(f)\neq 1$,
then my guess is trivially correct (in so far as the stated classification is trivial). However, what seems
not entirely trivial is that even more is true, in some sense. To flesh out this cryptic comment, we start by
recalling the situation in genus zero. Here, after some change of variables, one essentially reduces to
equations of the form
$$ax^2+by^2=c$$
were $a,b,c$ are positive. I am sure you are familiar then with the method of sweeping lines, whereby all the
solutions can be constructed from just one. On the other hand,  having an existence oracle will certainly
enable us to find the solution when it exists, just by exhaustively enumerating pairs until one is found. In
short, the existence oracle already provides us with a method for `constructing' all solutions.

A few years ago, I was happy to discover that a similar phenomenon occurs when the genus is at least two. That
is to say
\begin{observation}
For equations of genus at least two, relative to the existence problem, the full solution set is computable as
a function of $f$.
\end{observation}
Incidentally, the computability of the solution set for curves of higher genus is one of the two most
important questions regarding the arithmetic of curves, the other being the BSD conjecture. Usually,
number-theorists like to consider computability in a specific form, like a specific bound on the size of the
(numerators and denominators of the) solutions in terms of some simple algebraic invariants of $f$. This is
the subject of the effective Mordell conjecture, or Vojta's conjecture, or the ABC conjecture, and so on. Our
theorem proceeds in a different direction, and merely constructs an algorithm relative to the existence
oracle. In particular, no bounds on size, even relative ones, are obvious in this approach. My motivation for
pointing out this relative computability in fact was to emphasize a different perspective from the usual
arithmetic-geometric one, since an algorithm can be quite different from a special kind of formula. An
additional remark is that if a more powerful decision oracle is given,  for example, allowing equations and
inequations, then the theorem is trivial. The only interest comes from restricting the input to the forms
stated.

Let us briefly recall how one proves the theorem. Given an equation, assume the oracle tells us that a
solution exists. Now search in increasing order of size for a solution until we find one $(p,q)$. Then do a
quick additional search for more solutions with first coordinate $p$. Denote by $S$ the solution set we've found this
way. Now consider the two equations in three variables
$$f(x,y)=0, (x-p)z-1=0$$
It is easy to see that the solutions to this system is in bijection with the solutions of $f=0$ minus the ones
that we've found. In geometric language, we have embedded our curve into three space in such a way that the
set $S$ has been sent out to infinity. That is, if $X\subset \A^2$ is our curve, then we have embedded $X$
into $\P^3$ so that if $H$ denotes the plane at infinity then $X\cap H=S$. Now consider projections $\pi:
\P^3\setminus \{c\} \ra \P^2$ from some rational point $c\in H$. Then $H$ will map to a line in $\P^2$, which
we can then use as the new line at infinity. Furthermore, there always exists a  projection such that $X$ is
mapped birationally onto its image. But the tricky point is that we would like $X\setminus S$ to be mapped
{\em bijectively on rational points} onto an affine plane curve $X'$ which will then be defined by a new
equation $h(x,y)=0$ with exactly $|S|$ solutions less than $f=0$. Now the existence of a $c$ that will do the
trick is a  nice consequence of Hilbert's irreducibility theorem and rather pleasant plane geometry. Hilbert's
theorem comes in because it may not be possible to find a projection that is bijective on all points, but it
is always possible to find one bijective on rational points. Furthermore, for any given point $c$, it is
possible to check algorithmically whether or not the projection from $c$ satisfies this criterion. One need
only see if it lies on at most finitely many secants to $X$ and then, actually find those secants and check if
any of them are rational. All of this can be achieved using standard computational algebra programs. In this
manner, searching exhaustively, one locates the desired $c$.  Now one applies the oracle to $h=0$ and
proceeds.

The point of this discussion is that for curves of genus different from 1, the existence oracle is indeed very
powerful. Not only does it give an oracle for the finiteness problem, it provides us with a computable
function for the solution set, where the genus zero case of course is in a somewhat vague sense. By the way,
in the genus zero case, the existence problem is also decidable in a strong form.

Let us turn now to curves of genus one. As far as the decision problems themselves are concerned, this case is
the most interesting since the status of my guess is not obvious just from the classification. Therefore, I
was very pleased when Lou Van Den Dries pointed out to me that
\begin{thm}(Van Den Dries)
For curves of genus 1, the guess is correct.
\end{thm}
Combined then with the previous observations, we conclude that the guess is true for all curves.

The idea is that if the oracle tells us a solution exists, and we find one, then we are in the situation of an
elliptic curve. (Here, I ignore the subtlety that the curve may be singular. This eventuality can easily be
accommodated.) Then our equation has infinitely many solutions if and only if the elliptic curve has a
rational point of infinite order. But the rational points of finite order can be readily found using a theorem
of Nagell and Lutz. Having done that, we can eliminate all of them using the same trick as that outlined above
for higher genus curves. Then applying the oracle to the new equation finishes the job. In spite of the
simplicity of this remark, the result is that we have a rather tidy picture in the case of curves. Recall
again that the existence of points of infinite order on an elliptic curve is very hard to determine. An
algorithm for achieving this is a major consequence of the BSD conjecture. So it is rather interesting that an
existence oracle serves the same purpose. A few years ago, I asked John Tate his opinion on the existence
problem for genus one curves, whereupon he replied that it should be on the same order of difficulty as BSD.
In fact, it is a somewhat subtle fact (requiring the 2004 thesis of Alexander Perlis, a student at Arizona)
that BSD also does give us an existence oracle {\em and} and a finiteness oracle for arbitrary curves of genus
one. That is to say, we have the implications:

$$ \begin{array}{rcc}  & &\mbox{existence oracle} \\ & \nearrow & \\
\mbox{BSD} & & \\
& \searrow & \\
& & \mbox{finiteness oracle} \end{array}$$ which we have now filled into a triangle
$$ \begin{array}{rcc}  & &\mbox{existence oracle} \\ & \nearrow & \\
\mbox{BSD} & & \downarrow\\
& \searrow & \\
& & \mbox{finiteness oracle} \end{array}$$ The important point here is that the vertical arrow does not
require BSD. Since the finiteness oracle (at least for elliptic curves) is often thought of as {\em the} major
application of BSD, the arrow we have filled in can be taken as  vindication of Tate's intuition. In other
words, it is something of a weak implication $$\mbox{existence oracle} \cdots \rightarrow \mbox{BSD}.$$

I hope the above discussion has already given you some sense of why I feel the existence oracle to be a very
powerful thing in the study of rational solutions. But I should still explain a little bit the reasoning
behind my guess. It has to do with the geometric structure underlying the existence of infinitely many points.
There is a conjecture of Lang predicating that most of the rational points on any given variety are
concentrated inside a specific geometric locus. That is, given a variety $X$, denote by $E\subset X$ the
Zariski closure in $X$ of the images of all non-constant rational maps from group varieties. (Here, I mean
varieties that have the structure of an algebraic group over $\C$.) This is called by Lang the {\em
exceptional set} of $X$. Then Lang's conjecture says that $X\setminus E$ has finitely many rational points. If
true, one need only examine $E$ to decide the finiteness question. Of course, there is the question of
algorithmically finding $E$ from $X$, but the overall picture of the classification theory of algebraic
varieties makes it hard to believe that $E$ might not be a computable function of $X$. (Here, I am again
revealing my own naivete as an `ordinary' mathematician.) For $E$ itself, well, the structure of the decision
problem is not entirely clear. However, group varieties themselves do have the property that a decision oracle
gives a finiteness oracle, either for trivial reasons (they might be geometrically rational varieties) or for
deep ones (BSD). The precise location of my guess then is that this property should be inherited by images of
group varieties and eventually the whole exceptional set. At some later time, I hope to think about this issue
seriously enough to remove the tentative nature of this discussion. In the meanwhile, I hope it is at least
clear enough that the relative decidability question for rational points really does tie in to central
problems of Diophantine geometry.

Perhaps I can close with one other question that I'm sure some of you have asked already in the course of this
lecture, if not before. What about the converse implication? That is, either for integral or rational
solutions, is the existence problem decidable relative to the finiteness problem? Because of the dummy variable trick, there is of course a positive answer for
the general problem. But it is of
interest to ask this question in a limited context, for example,
after fixing the number of variables. Even then, I would guess that
this question is less natural
than the one we've been discussing from the perspective of recursion theory. 
However, 
for the arithmetic of rational points on curves, a positive answer would be extremely powerful. After all, we {\em do} have a finiteness oracle for curves of genus not equal to 1,
and for all curves if we assume BSD.

\end{document}